\input ./style/arxiv-general.cfg
\documentclass[aos,MSNbibl,nameyear,dvips]{arximspdf}
\makeatletter
   \@ifpackageloaded{graphicx}{}{\usepackage{graphicx}}
\makeatother
\doi{10.1214/15-AOS1317}% Updated by VTEXPTS2LaTeX.exe, 02.03.2015
\volume{43}
\issue{4}
\pubyear{2015}
\firstpage{1596}
\lastpage{1616}
\docsubty{FLA}

\makeatletter
\newcommand{\tr}{\operatorname{tr}}

\newtheorem{lemma}{Lemma}
\newproclaim{ex}{Example}
\newtheorem{theorem}{Theorem}
\newproclaim{remark}{Remark}
\newproclaim{definition}{Definition}

\newproclaim{example}{Example}
\newtheorem{proposition}{Proposition}

\makeatother

\begin{document}
\begin{frontmatter}

\title{Optimal designs for the proportional interference~model}
\runtitle{Optimal designs for the proportional interference model}

\begin{aug}
% Corresponding author: wei zheng - weizheng@iupui.edu% Updated by
%VTEXPTS2LaTeX.exe, 02.03.2015 10:14
\author[A]{\fnms{Kang}~\snm{Li}\thanksref{m1}\ead[label=e1]{likang1115@163.com}},
\author[B]{\fnms{Wei}~\snm{Zheng}\corref{}\thanksref{m2}\ead[label=e2]{weizheng@iupui.edu}}
\and
\author[A]{\fnms{Mingyao}~\snm{Ai}\thanksref{m1,T1}\ead[label=e3]{myai@pku.edu.cn}}
\runauthor{K. Li, W. Zheng and M. Ai}
\affiliation{Peking University\thanksmark{m1} and Indiana
University-Purdue University Indianapolis\thanksmark{m2}}
%\dedicated{}
\address[A]{LMAM, School of Mathematical Sciences \\
\quad and Center for Statisitcal Science\\
Peking University\\
Beijing 100871\\
China\\
\printead{e1}\\
\phantom{E-mail:\ }\printead*{e3}}
\address[B]{Department of Mathematical Science\\
Indiana University-Purdue University\\
Indianapolis, Indiana 46202\\
USA\\
\printead{e2}}
\end{aug}
\thankstext{T1}{Supported by NSFC Grants 11271032 and 11331011, BCMIIS
and LMEQF.}

% HISTORY:
%
\received{\smonth{10} \syear{2014}}% Updated by VTEXPTS2LaTeX.exe,
%02.03.2015 10:14
%
\revised{\smonth{1} \syear{2015}}% Updated by VTEXPTS2LaTeX.exe,
%02.03.2015 10:14

% ABSTRACT
%
\begin{abstract}
The interference model has been widely used and studied in block
experiments where the treatment for a particular plot has effects on its
neighbor plots. In this paper, we study optimal circular designs for
the proportional interference model, in which the neighbor effects of a
treatment are proportional to its direct effect. Kiefer's equivalence
theorems for estimating both the direct and total treatment effects are
developed with respect to the criteria of A, D, E and T. Parallel
studies are carried out for the undirectional model, where the neighbor
effects do not depend on whether they are from the left or right.
Moreover, the connection between optimal designs for the directional and
undiretional models is built. Importantly, one can easily develop a
computer program for finding optimal designs based on these theorems.
\end{abstract}

% KEYWORDS
\begin{keyword}[class=AMS]
\kwd[Primary ]{62K05}
%\kwd{}
\kwd[; secondary ]{62J05}
\end{keyword}
\begin{keyword}
\kwd{Approximate design theory}
\kwd{equivalence theorem}
\kwd{interference model}
\kwd{optimal design}
\kwd{proportional model}
\end{keyword}
\end{frontmatter}

\section{Introduction}

In many agricultural experiments, the treatment assigned to a
particular plot could also have effects on its neighbor plots. This is well
recognized in literature. See \citet{r11}, \citet{r22}, \citet{r7}, \citet{r27}, \citet{r19} and \citet{r20}, for examples.
To adjust the biases caused by
these neighbor effects, the interference model is widely adopted. In a
block design with $n$ blocks of size $k$ and $t$ treatments, the
response, $y_{dij}$, observed from the $j$th plot of block $i$ is decomposed
into the following items:
\begin{eqnarray}
\label{eqn729}
y_{dij}&=&\mu+\beta_i+\tau_{d(i,j)}+
\gamma_{d(i,j-1)}+\rho _{d(i,j+1)}+\varepsilon_{ij},
\end{eqnarray}
where the subscript $d(i,j)$ denotes the treatment assigned to the
$j$th plot of block $i$ by the design $d\dvtx \{1,2,\ldots,n\}\times
\{1,2,\ldots,k\}\rightarrow\{1,2,\ldots,t\}$. Furthermore, $\mu$ is the
general mean, $\beta_i$ is the effect of block $i$, $\tau_{d(i,j)}$
is the
direct effect of treatment $d(i,j)$, $\lambda_{d(i,j-1)}$ is the
neighbor effect of treatment $d(i,j-1)$ from the left and $\rho
_{d(i,j+1)}$ is
the neighbor effect of treatment $d(i,j+1)$ from the right. At last,
$\varepsilon_{ij}$ is the error term.

\citet{r24} studied optimal designs under model (\ref
{eqn729}) for estimating the direct treatment effect when $k=3$ or
$4$. The
latter was extended to $5\leq k\leq t$ by \citet{r25}.
\citet{r33} recently provided a unified framework in deriving optimal
designs for general values of $k$ and $t$, with an arbitrary structure
of the within-block covariance matrix. On the other hand, \citet{r4} studied the optimal designs under the same model,
however, for estimating the total treatment effect which is the
summation of
the direct and neighbor effects. This line of research was extended by
\citet{r1}, \citet{r2} and \citet{r14}. See also \citet{r19}, \citet{r12}, Filipiak and Markiewicz
(\citeyear{r16}, \citeyear{r17}, \citeyear{r18}) and \citet{r15} among others for relevant works on
optimal designs.

In this paper, we shall consider the proportional interference model,
where the neighbor effects are proportional to the direct treatment
effect, that is, $\gamma_i=\lambda_1\tau_i$ and $\rho_i=\lambda
_2\tau
_i$, $1\leq i\leq t$, for unknown constants $\lambda_1$ and $\lambda
_2$. This
is reasonable for many applications since an effective treatment
typically has large impacts on its neighbor plots. In fact, \citet{r11} has proposed such model with $\lambda_1=\lambda_2$. A model with
this restriction is said to be undirectional; otherwise, it is
directional. Yet, there is no literature on optimal designs under
either of these two models according to the best knowledge of the authors.
Meanwhile, optimal crossover designs under a similar proportional model
have been studied by \citet{r23}, \citet{r6}, \citet{r9} and \citet{r31}. By their
enlightenment, the nonlinear terms $\lambda_1\tau_i$ and $\lambda
_2\tau
_i$ in the
proportional interference model can be handled in the same fashion. We
are interested in finding the optimal designs for estimating the direct
and total treatment effects, respectively, under either of the
directional and undirectional models.

Let $Y_d$ be the vector of responses organized block by block. Now we
can write the proportional interference model as follows:
\begin{eqnarray}
\label{eqn728}
Y_d&=&1_{nk}\mu+U\beta+(T_d+
\lambda_1L_d+\lambda_2R_d)\tau
+\varepsilon,
\end{eqnarray}
where $1_{nk}$ represents a vector of ones with length $nk$, $\beta
=(\beta_1,\ldots,\beta_n)'$, $\tau=(\tau_1,\ldots,\tau_t)'$ and
$U=I_n\otimes
1_k$. Here, $I_n$ is the identity matrix of order $n$; $\otimes$
denotes the Kronecker product and $'$ means the transposition. Also, $T_d$,
$L_d$ and $R_d$ represent the design matrices for the direct, left
neighbor and right neighbor effects, respectively. Throughout the
paper, we
consider circular designs, for which $d(i,0)=d(i,k)$ and
$d(i,k+1)=d(i,1)$, $1\le i\le n$. Hence we have $L_d=(I_n\otimes H)
T_d$ and
$R_d=(I_n\otimes H') T_d$, where $H=(\mathbb{I}_{i=j+1 (\operatorname{mod} \, \, k)})_{1\leq
i,j\leq k}$ with the indicator function $\mathbb{I}$. For the random
error term
$\varepsilon$, we assume that ${\mathbb{E}}(\varepsilon)=0$ and
$\operatorname{Var}(\varepsilon
)=I_n\otimes\Sigma$, where $\Sigma$ is an arbitrary $k\times k$
positive definite within-block covariance matrix.

The rest of the paper is organized as follows. Sections~\ref{secOptimaldesignsfordirecteffect} and \ref{secOptimaldesignsfortotaleffect} investigate the optimal designs for estimating the direct and
total treatment effects, respectively, under the proportional interference
model.
Kiefer's equivalence theorems are given with respect to A, D, E and T
criteria therein.
Section~\ref{secUndirectionalsideeffects} carries out parallel
studies for
the undirectional model. Moreover, the connection between optimal
designs for the two models is built. Section~\ref{secexamples} illustrates
these theorems through several examples. Section~\ref{secconclusion}
concludes the paper with some discussions.

\section{Optimal designs for direct treatment effect}\label{secOptimaldesignsfordirecteffect}
For any matrix $G$, define $G^-$ as a generalized inverse of $G$ and
the projection operator $\operatorname{pr}^{\bot}G=I-G(G'G)^-G'$. Let
$\tilde
{U}=(I_n\otimes
\Sigma^{-1/2})U$, $\tilde{T}_d=(I_n\otimes\Sigma^{-1/2})T_d$,
$\tilde
{L}_d=(I_n\otimes\Sigma^{-1/2})L_d$ and $\tilde{R}_d=(I_n\otimes
\Sigma^{-1/2})R_d$. The Fisher information matrix for the direct
treatment effect $\tau$ under model (\ref{eqn728}) is
\begin{eqnarray*}
C_d(\tau)&=&(\tilde{T}_d+\lambda_1
\tilde{L}_d+\lambda_2\tilde {R}_d)'
\operatorname{pr}^{\bot} (\tilde{U}|\tilde{L}_d\tau|
\tilde{R}_d\tau) (\tilde{T}_d+\lambda _1
\tilde {L}_d+\lambda_2\tilde{R}_d).
\end{eqnarray*}
For notational convenience, let $M_{x,y,z}=x\tilde{T}_d+y\tilde
{L}_d+z\tilde{R}_d$ for any values of $x,y$ and $z$. By setting
$\lambda
_0=1$, we
have
\begin{eqnarray}
C_d(\tau)&=&M_{1,\lambda_1,\lambda_2}'\operatorname{pr}^{\bot}(
\tilde {U}|\tilde{L}_d\tau |\tilde{R}_d
\tau)M_{1,\lambda_1,\lambda_2}
\nonumber
\\
&=&M_{1,\lambda_1,\lambda_2}'\operatorname{pr}^{\bot}(\tilde
{U})M_{1,\lambda_1,\lambda
_2}-M_{1,\lambda_1,\lambda_2}'\operatorname{pr}^{\bot}(
\tilde{U}) (\tilde {L}_d\tau|\tilde {R}_d\tau)
\nonumber
\\
\label{eqn804}
&&\times\bigl[(\tilde{L}_d\tau|\tilde{R}_d
\tau)'\operatorname{pr}^{\bot
}(\tilde{U}) (\tilde {L}_d
\tau|\tilde{R}_d\tau)\bigr]^-(\tilde{L}_d\tau|
\tilde{R}_d\tau )'\operatorname{pr}^{\bot} (
\tilde{U})M_{1,\lambda_1,\lambda_2}
\\
&=&\sum^2_{i=0}\sum
^2_{j=0}\lambda_i\lambda_jC_{dij}-A_d'
\bigl(\tau 'C_{dij}\tau\bigr)_{1\leq i,j\leq2}^-A_d,
\nonumber
\\
\nonumber
A_d&=& \Biggl(\sum_{i=0}^2
\lambda_iC_{di1}\tau\bigg|\sum_{i=0}^2
\lambda _iC_{di2}\tau \Biggr)',
\nonumber
\end{eqnarray}
where $C_{dij}=G_i'(I_n\otimes\tilde{B})G_j$, $0\leq i,j\leq2$, with
$G_0=T_d$, $G_1=L_d$, $G_2=R_d$ and
$\tilde{B}=\Sigma^{-1}-\Sigma^{-1}1_k1_k'\Sigma^{-1}/1_k'\Sigma
^{-1}1_k$. In particular, if $\Sigma$ is a matrix of {\it type-H},
that is,
$\Sigma=I_k+\eta1_k'+1_k\eta'$ with a vector $\eta$ of length $k$, we
have $\tilde{B}=\operatorname{pr}^{\bot}(1_k):=B_k$ [\citet{r26}].
Examples of type-H
matrices include the identity matrices and completely symmetric matrices.

One major objective of design theorists is to find a design with
maximum information matrix. Following Kiefer (\citeyear{K75}), we shall try to find
the designs which maximize $\Phi(C_d(\tau))$, where $\Phi$ satisfies
the following three conditions:
\begin{longlist}[(C.3)]
\item[(C.1)] $\Phi$ is concave;
\item[(C.2)] $\Phi(S'CS)=\Phi(C)$ for any permutation matrix $S$;
\item[(C.3)] $\Phi(bC)$ is nondecreasing in the scalar $b>0$.
\end{longlist}

Note that $C_d(\tau)$ depends on the true value of $\tau$ itself, and
thus the choice of optimal designs. Following \citet{r23}, \citet{r6} and \citet{r31}, we adopt the
Bayesian type criterion
\begin{eqnarray}\label{eqn0824}
\phi_g(d)&=&\int\Phi\bigl(C_d(\tau)\bigr)g(\tau)
\,\mbox{d}(\tau) =\mathbb {E}_g\bigl(\Phi \bigl(C_d(\tau)
\bigr)\bigr),
\end{eqnarray}
where $g$ is the prior distribution of $\tau$ and is assumed to be
exchangeable throughout the paper. A design is said to be optimal if
and only
if it achieves the maximum of $\phi_g(d)$ among all designs for given
$g$, $\Phi$, $\lambda_1$ and $\lambda_2$. Furthermore, if a design
maximizes $\phi_g(d)$ for any $\Phi$, it is also said to be {\it
universally optimal}.

In this paper we consider four popular criteria for finding optimal
designs. For a $t\times t$ matrix $C$ with eigenvalues $0=a_0\leq
a_1\leq
a_2\leq\cdots\leq a_{t-1}$, define the criterion functions as
\begin{eqnarray*}
\Phi_A(C)&=&(t-1) \Biggl(\sum^{t-1}_{i=1}a_i^{-1}
\Biggr)^{-1},
\\
\Phi_D(C)&=& \Biggl(\prod_{i=1}^{t-1}a_i
\Biggr)^{1/(t-1)},
\\
\Phi_E(C)&=&a_1,
\\
\Phi_T(C)&=&(t-1)^{-1}\sum^{t-1}_{i=1}a_i.
\end{eqnarray*}
A design is said to be ${\mathcal A}_g$-optimal if it maximizes $\phi
_g(d)$ with $\Phi=\Phi_A$ in (\ref{eqn0824}). The ${\mathcal D}_g$-,
${\mathcal E}_g$-
and ${\mathcal T}_g$-optimality of a design are similarly defined.

Let $\Omega_{n,k,t}$ denote the set of all possible block designs with
$n$ blocks of size $k$ and $t$ treatments. A design in $\Omega_{n,k,t}$
could be considered as a result of selecting $n$ elements from the set,
${\mathcal S}$, of all possible $t^k$ block sequences with replacement. For
each $s\in{\mathcal S}$, we define the sequence proportion
$p_s=n_s/n$, where $n_s$ is the number of replications of $s$ in the
design. For given
$n$, a design is determined by the measure $\xi=(p_s,s\in{\mathcal
S})$. If $p_s>0$, then $s$ is a supporting sequence of~$\xi$. In approximate
design theory, we search for optimal measures in the space of
${\mathcal P}=\{ \xi\dvtx  \sum_{s\in{\mathcal S}} p_s=1,p_s\geq0\}$. If
such a
measure happens to fall within the subset ${\mathcal P}_n=\{\xi\in
{\mathcal P}\dvtx  n\xi\mbox{ is a vector of integers}\}$, then we derive
an exact
design which is optimal among $\Omega_{n,k,t}$.

Let $C_{sij}$ be the matrix $C_{dij}$ when the design $d$ is
degenerated to a single sequence $s$ for $0\leq i,j\leq2$. Then we have
$C_{dij}=n\sum_{s\in{\mathcal S}} p_sC_{sij}$. By equation (\ref
{eqn804}), we have
\begin{eqnarray}
C_{d}(\tau)&=&nC_{\xi}(\tau),
\nonumber
\\
C_{\xi}(\tau)&=&\sum^2_{i=0}\sum
^2_{j=0}\lambda_i
\lambda_jC_{\xi
ij}-A_{\xi}'\bigl(
\tau'C_{\xi ij}\tau\bigr)_{1\leq i,j\leq2}^-A_{\xi},
\nonumber
\\[-8pt]
\label{eqn2014805}
\\[-8pt]
\nonumber
C_{\xi ij}&=&\sum_{s\in{\mathcal S}} p_sC_{sij},
\nonumber
\\
A_\xi&=& \Biggl(\sum_{i=0}^2
\lambda_iC_{\xi i1}\tau\bigg|\sum_{i=0}^2
\lambda_iC_{\xi i2}\tau \Biggr)'.
\nonumber
\end{eqnarray}
Note that $C_{\xi}(\tau)$ is independent of $n$. Here we call $C_{\xi
}(\tau)$ the information matrix of the measure $\xi$. Furthermore, by
noting that the four criterion functions satisfy $\Phi(nC)=n\Phi(C)$,
we have
\begin{eqnarray}
\phi_g(d)&=&n\phi_g(\xi),
\nonumber
\\[-8pt]
\label{eqn605}
\\[-8pt]
\nonumber
\phi_g(\xi)&=&\int\Phi\bigl(C_{\xi}(\tau)\bigr)g(\tau)
\,\mbox{d}(\tau).
\nonumber
\end{eqnarray}
Equation (\ref{eqn605}) indicates that the number of blocks $n$ is
irrelevant to the search of approximate optimal designs. In the sequel,
we shall focus on finding the optimal measures which maximize $\phi
_g(\xi)$ among ${\mathcal P}$.

Let $\mathcal O$ denote the set of all $t!$ permutation operators on $\{
1,2,\ldots,t\}$. For any $\sigma\in\mathcal O$ and $s=(t_1,\ldots
,t_k)$ with
$1\le t_i\le t$, define $\sigma s=(\sigma(t_1),\ldots,\sigma(t_k))$. A
measure is said to be {\it symmetric} if it is invariant under treatment
relabeling, that is, $\sigma\xi=\xi$ for all $\sigma\in\mathcal O$,
where $\sigma\xi=(p_{\sigma^{-1}s},s\in{\mathcal S})$. By adopting similar
arguments to Corollary 1 in \citet{r31}, we get the following result.

\begin{proposition}\label{prop730}
In approximate design theory, given any values of $\lambda_1$ and~$\lambda_2$, and the exchangeable prior distribution $g$ of $\tau$,
for any
measure $\xi$ there exists a symmetric measure, say $\xi^*$, such that
\begin{eqnarray*}
\phi_{g}(\xi)&\leq&\phi_{g}\bigl(\xi^*\bigr).
\end{eqnarray*}
\end{proposition}

Proposition \ref{prop730} indicates that an optimal measure in the
subclass of symmetric measures is automatically optimal among
${\mathcal P}$. The merit of such a result is
that the form of the information matrix for a symmetric measure is
usually feasible to be calculated explicitly. In fact there is a larger
subclass of measures with the same convenience. We say a measure is
{\it pseudo symmetric} if $C_{\xi ij}$, $0\le i,j\le2$ are all
completely symmetric. A symmetric measure is also pseudo symmetric
[\citet{r26}]. It is easy to verify that the column and row sums of
$C_{\xi ij}$'s are all zero. Hence, for any pseudo symmetric measure we
have $C_{\xi ij}=c_{\xi ij}B_t/(t-1)$, $0\leq i,j\leq2$, where $c_{\xi
ij}=\tr(C_{\xi ij})$. Now let $\ell=(1,\lambda_1,\lambda_2)'$,
$V_{\xi
}=(c_{\xi ij})_{0\leq i,j\leq2}$, $Q_{\xi}=(c_{\xi ij})_{1\leq
i,j\leq
2}$ and
\begin{eqnarray}
\label{eqn1109}
q_\xi^*&=&c_{\xi00}-\pmatrix{
c_{\xi01} & c_{\xi02}}
  Q_{\xi}^- \pmatrix{
c_{\xi10}
\cr
c_{\xi20}
}.
\end{eqnarray}

\begin{proposition}\label{prop7302}
For a pseudo symmetric measure $\xi$, the information matrix $C_{\xi
}(\tau)$ has eigenvalues of $0$, $(t-1)^{-1}q_{\xi}^*$ and
$(t-1)^{-1}\ell'V_{\xi}\ell$ with multiplicities of~$1$, $1$ and $t-2$,
respectively. Moreover we have $q_{\xi}^*\leq\ell'V_{\xi}\ell$.
\end{proposition}
\begin{pf}
Due to $1_t'\tau=0$, we have $B_t\tau=\tau$ and $\tau'C_{\xi
ij}\tau
=c_{\xi ij}\tau'\tau/(t-1)$. In view of (\ref{eqn2014805}), we obtain
\begin{eqnarray*}
(t-1)C_{\xi}(\tau)&=&\sum_{i=0}^2
\sum_{j=0}^2\lambda_i\lambda
_jc_{\xi
ij}B_t-a\bigl(\tau'\tau
\bigr)^{-1}\tau\tau',
\\
a&=& \Biggl(\sum_{i=0}^2
\lambda_ic_{\xi i1}\bigg|\sum_{i=0}^2
\lambda _ic_{\xi i2} \Biggr) Q_{\xi}^- \pmatrix{
\displaystyle\sum_{j=0}^2
\lambda_jc_{\xi1j}
\vspace*{3pt}\cr
\displaystyle\sum_{j=0}^2\lambda_jc_{\xi2j}
}.
\end{eqnarray*}

Let $\{x_1,\ldots,x_{t-2}\}$ be the orthogonal basis that is orthogonal to
both $1_t$ and $\tau$. Then $\{1_t,\tau,x_1,\ldots,x_{t-2}\}$ forms the
eigenvectors of $C_{\xi}(\tau)$. The corresponding eigenvalues are $0$,
$(t-1)^{-1}(\sum_{i=0}^2\sum_{j=0}^2\lambda_i\lambda_jc_{\xi
ij}-a)$ and
$(t-1)^{-1}\sum_{i=0}^2\sum_{j=0}^2\lambda_i\lambda_jc_{\xi ij}$ with
multiplicities of $1$, $1$ and $t-2$, respectively. The proof is concluded
in view of
\begin{eqnarray*}
\ell'V_{\xi}\ell&=&\sum_{i=0}^2
\sum_{j=0}^2\lambda_i\lambda
_jc_{\xi
ij},
\\
a&=&\ell' \pmatrix{
c_{\xi01} &
c_{\xi02}
\cr
c_{\xi11} & c_{\xi12}
\cr
c_{\xi21} & c_{\xi22}
}
  Q_{\xi}^- \pmatrix{
c_{\xi01} & c_{\xi11} & c_{\xi21}
\cr
c_{\xi02} & c_{\xi12} & c_{\xi22}
}
  \ell\geq0
\end{eqnarray*}
and (\ref{eqn1109}).
\end{pf}

By Proposition \ref{prop7302}, it is seen that $\phi_g(\xi)=\Phi
(C_{\xi
}(\tau))$ for any pseudo symmetric measure under the four criterion
functions. Hence $g$ is irrelevant to the determination of optimal
pseudo symmetric measures for the four criteria.

\begin{lemma}\label{lemma1109}
Except for measures with each supporting sequence consisting of only
one treatment, we have $c_{\xi ii}>0$ for $i=0,1,2$. If $\det(Q_{\xi})>0$,
then $q_{\xi}^*=\det(V_{\xi})/\det(Q_{\xi})$, where $\det(\cdot
)$ means
the determinant of a matrix. Otherwise, $q_{\xi}^*=c_{\xi00}-c_{\xi
01}^2/c_{\xi11}=c_{\xi00}-c_{\xi02}^2/c_{\xi22}$.
\end{lemma}
\begin{pf}
%Since $\proj U=I_n\otimes B_k$ and $HB_kH'=H'B_kH=B_k$, we have $c_{
%\xi00}=c_{\xi11}=c_{\xi22}$. Note that $c_{\xi00}=0$ is equivalent
%to
%$\proj(U)T_d=0$, and the latter is only possible when each block
%repeats the same treatment throughout the $k$ plots.
%For a non-circular design, we know that $\proj(U)L_d1_t\neq0$ and $
%\proj(U)R_d1_t\neq0$, and hence $c_{d11}c_{d22}>0$.
%The rest of the lemma follows by straight forward calculations.
Note that $\tilde{B}$ is nonnegative definite. So $c_{\xi00}=\sum_{s\in
\mathcal S}p_s \tr(T_s'\tilde{B}T_s)$ $\ge0$, where $T_s$ is the
matrix $T_d$
when $d$ is degenerated to a single sequence $s$. If $c_{\xi00}=0$, we
have $T_s' \tilde{B}T_s=0$ and thus $\tilde{B}T_s=0$ for any supporting
sequence $s$. It is known that $\tilde{B}x=0$ if and only if $x$ is a
multiple of $1_k$ [\citet{r26}]. This is only possible when each
supporting sequence repeats the same treatment throughout the $k$
plots. For $c_{\xi11}$ and $c_{\xi22}$, we have the similar
arguments. The
rest of the lemma follows by straightforward calculations.
\end{pf}

From the proof of Lemma \ref{lemma1109}, we know $V_{\xi}=0$ if each
supporting sequence of $\xi$ consists of only one treatment. There is no
information gathered from such measures regarding $\tau$, and hence it
is impossible to be optimal. In the subsequent arguments, we neglect such
measures by default.

Define the quadratic functions $q_{\xi}(x)=c_{\xi00}+2c_{\xi
01}x+c_{\xi11}x^2$ and $x_{\xi}=-c_{\xi01}/c_{\xi11}$. Then we have
$q_{\xi}(x_{\xi})=c_{\xi00}-c_{\xi01}^2/c_{\xi11}$. Let $c_{s
ij}=\tr(C_{s ij})$, $V_{s}=(c_{s ij})_{0\leq i,j\leq2}$, $Q_{s}=(c_{s
ij})_{1\leq i,j\leq2}$ and\vspace*{1pt} $q_s(x)=c_{s00}+2c_{s01}x+c_{s11}x^2$.
Clearly, $c_{\xi ij}=\sum_{s\in\mathcal S} p_s c_s$, $V_{\xi}=\sum_{s\in\mathcal S}
p_s V_s$, $Q_{\xi}=\sum_{s\in\mathcal S} p_s Q_s$ and $q_{\xi
}(x)=\sum_{s\in\mathcal S} p_s q_s(x)$.

\begin{theorem}\label{thm1110}
In estimating $\tau$ under model (\ref{eqn728}), a pseudo symmetric
measure $\xi$ is optimal in the following cases. In each case, the
$p_s$ in
$\xi$ is positive only if $s$ reaches the maximum therein.
\begin{longlist}[(iii)]
\item[(i)] If $\det(Q_{\xi})=0$, then $\xi$ is ${\mathcal A}_g$-optimal if
and only if
\begin{eqnarray*}
&& \max_{s\in{\mathcal
S}}\frac{q_{\xi}(x_{\xi})^{-2}q_s(x_{\xi})+(t-2)(\ell'V_{\xi}\ell
)^{-2}\ell'V_s\ell}{q_{\xi}(x_{\xi})^{-1}+(t-2)(\ell'V_{\xi}\ell
)^{-1}}=1.
\end{eqnarray*}
If $\det(V_{\xi})>0$, then $\xi$ is ${\mathcal A}_g$-optimal if and
only if
\begin{eqnarray*}
&& \max_{s\in{\mathcal S}}\frac{r_s
\det(Q_{\xi})/\det(V_{\xi})+(t-2)(\ell'V_{\xi}\ell)^{-2}\ell
'V_s\ell
}{\det(Q_{\xi})/\det(V_{\xi})+(t-2)(\ell'V_{\xi}\ell)^{-1}}=1.
\end{eqnarray*}
\item[(ii)]  If $\det(Q_{\xi})=0$, then $\xi$ is ${\mathcal D}_g$-optimal if
and only if
\begin{eqnarray}
\label{eqn11098}
\max_{s\in{\mathcal S}} \biggl(\frac{1}{t-1}
\frac{q_s(x_{\xi
})}{q_{\xi
}(x_{\xi})}+\frac{t-2}{t-1}\frac{\ell'V_s\ell}{\ell'V_{\xi}\ell
} \biggr)&=&1.
\end{eqnarray}
If $\det(V_{\xi})>0$, then $\xi$ is ${\mathcal D}_g$-optimal if and
only if
\begin{eqnarray}
\label{eqn8016} \max_{s\in{\mathcal S}} \biggl(\frac{r_s}{t-1}+
\frac{t-2}{t-1}\frac
{\ell
'V_s\ell}{\ell'V_{\xi}\ell} \biggr)&=&1.
\end{eqnarray}
\item[(iii)]  If $\det(Q_{\xi})=0$, then $\xi$ is ${\mathcal E}_g$-optimal if
and only if
\begin{eqnarray*}
\max_{s\in{\mathcal S}}\frac{q_s(x_{\xi})}{q_{\xi}(x_{\xi})}&=&1.
\end{eqnarray*}
If $\det(V_{\xi})>0$, then $\xi$ is ${\mathcal E}_g$-optimal if and
only if
\begin{eqnarray*}
\max_{s\in{\mathcal S}}r_s&=&1.
\end{eqnarray*}
\item[(iv)]  If $\det(Q_{\xi})=0$, then $\xi$ is ${\mathcal T}_g$-optimal if
and only if
\begin{eqnarray*}
\max_{s\in{\mathcal S}}\frac{q_s(x_{\xi})+(t-2)\ell'V_s\ell
}{q_{\xi
}(x_{\xi})+(t-2)\ell'V_{\xi}\ell}&=&1.
\end{eqnarray*}
If $\det(V_{\xi})>0$, then $\xi$ is ${\mathcal T}_g$-optimal if and
only if
\begin{eqnarray*}
\max_{s\in{\mathcal S}}\frac{r_s \det(V_{\xi})/\det(Q_{\xi
})+(t-2)\ell
'V_s\ell}{\det(V_{\xi})/\det(Q_{\xi})+(t-2)\ell'V_{\xi}\ell}&=&1.
\end{eqnarray*}
Here $r_s=\tr(V_sV_{\xi}^{-1})-\tr(Q_sQ_{\xi}^{-1})$.
\end{longlist}
\end{theorem}
\begin{pf}
Here we give only the proof for (ii) and the other three cases follow
similarly. First we would like to show that
\begin{eqnarray}
\label{eqn11092}
\det(V_{\xi})/\det(Q_{\xi})&
\leq&c_{\xi00}-c_{\xi01}^2/c_{\xi11}
\end{eqnarray}
whenever $\det(Q_{\xi})>0$. To see this, consider the following inequality:
\begin{eqnarray}
\label{eqn11093}
\pmatrix{
c_{\xi00} &
c_{\xi01}
\cr
c_{\xi10} & c_{\xi11}
} -\frac{1}{c_{\xi22}}\pmatrix{
c_{\xi02}
\cr
c_{\xi12}}
\pmatrix{
c_{\xi20} & c_{\xi21}
} \leq \pmatrix{
c_{\xi00} & c_{\xi01}
\cr
c_{\xi10} & c_{\xi11}
}.
\end{eqnarray}
The left (resp.,  right) hand side of (\ref{eqn11092}) is the Schur
complement of the left (resp., right) hand side of (\ref{eqn11093}), and
hence (\ref{eqn11092}) follows by the nondecreasing property of Schur
complement.

By the definition of ${\mathcal D}_g$-optimality, Propositions \ref
{prop730} and \ref{prop7302}, Lemma \ref{lemma1109} and inequality
(\ref{eqn11092}), a pseudo symmetric measure $\xi$ with $\det
(Q_{\xi
})=0$ is ${\mathcal D}_g$-optimal if and only if
\begin{eqnarray}
\label{eqn11094}
\lim_{\delta\rightarrow0} \frac{\psi[(1-\delta)\xi+\delta\xi
_0]-\psi
(\xi)}{\delta}&\leq&0
\end{eqnarray}
for any measure $\xi_0$, where $\psi(\xi)=\log(q_{\xi}(x_{\xi
}))+(t-2)\log(\ell'V_{\xi}\ell)$. Here we used the fact that
$V_{\xi_0^*}=V_{\xi_0}$ and hence $\psi(\xi_0^*)=\psi(\xi_0)$, where
$\xi_0^*$ is a symmetric measure defined by $\xi_0^*=\sum_{\sigma
\in
\mathcal
O}\sigma\xi_0/t!$. Direct calculations show that (\ref{eqn11094}) is
equivalent to
\begin{eqnarray}
\label{eqn11095}
\frac{1}{t-1}\frac{q_{\xi_0}(x_{\xi})}{q_{\xi}(x_{\xi})}+\frac
{t-2}{t-1}
\frac{\ell'V_{\xi_0}\ell}{\ell'V_{\xi}\ell}&\leq&1.
\end{eqnarray}
By reducing $\xi_0$ to a degenerate measure which puts all mass on a
single sequence~$s$, we have
\begin{eqnarray}
\label{eqn11096} \max_{s\in{\mathcal S}} \biggl(\frac{1}{t-1}
\frac{q_s(x_{\xi
})}{q_{\xi
}(x_{\xi})}+\frac{t-2}{t-1}\frac{\ell'V_s\ell}{\ell'V_{\xi}\ell
} \biggr)&\leq &1.
\end{eqnarray}
By letting $\xi_0=\xi$, we have equality in (\ref{eqn11095}) and hence
\begin{eqnarray}
\label{eqn11097} \max_{s\in{\mathcal S}} \biggl(\frac{1}{t-1}
\frac{q_s(x_{\xi
})}{q_{\xi
}(x_{\xi})}+\frac{t-2}{t-1}\frac{\ell'V_s\ell}{\ell'V_{\xi}\ell
} \biggr)&\geq&1
\end{eqnarray}
in view of $q_{\xi}(x)=\sum_{s\in{\mathcal S}}p_sq_s(x)$. Combining
(\ref{eqn11096}) and (\ref{eqn11097}), we obtain (\ref{eqn11098}).

%Note that the function $\log(q_{\xi}(x_{\xi}))+(t-2)\log(\ell'V_{\xi}
%\ell)$ is concave in $\xi$ globally.

For a pseudo symmetric measure $\xi$ with $\det(V_{\xi})>0$ and any
measure $\xi_0$, by the continuity of $\det(Q_{(1-\delta)\xi+\delta
\xi_0})$
in $\delta$, there exists a constant $\epsilon>0$ such that $\det
(Q_{(1-\delta)\xi+\delta\xi_0})>0$ for all $\delta\in(-\epsilon,
\epsilon)$.
Hence $\xi$ is ${\mathcal D}_g$-optimal if and only if
\begin{eqnarray}
\label{eqn8014} \lim_{\delta\rightarrow0} \frac{\varphi[(1-\delta)\xi+\delta\xi
_0]-\varphi(\xi)}{\delta}&\leq&0,
\end{eqnarray}
where $\varphi(\xi)=\log(\det(V_{\xi})/\det(Q_{\xi}))+(t-2)\log
(\ell
'V_{\xi}\ell)$. It is well known that
\begin{eqnarray}
\label{eqn8013}
\lim_{\delta\rightarrow0} \frac{\log(\det(V_{(1-\delta)\xi
+\delta\xi
_0}))-\log(\det(V_{\xi}))}{\delta}&=&\tr
\bigl(V_{\xi_0}V_{\xi}^{-1}\bigr)-3.
\end{eqnarray}
The same result holds for $Q_{\xi}$ except that the number $3$ in
(\ref
{eqn8013}) is replaced with $2$. By applying (\ref{eqn8013}) to
(\ref{eqn8014}) we have
\begin{eqnarray}
\label{eqn8015} \frac{\tr(V_{\xi_0}V_{\xi}^{-1})-\tr(Q_{\xi_0}Q_{\xi
}^{-1})}{t-1}+\frac
{t-2}{t-1}\frac{\ell'V_{\xi_0}\ell}{\ell'V_{\xi}\ell}&\leq&1.
\end{eqnarray}
Hence, for single sequences we have
\begin{eqnarray*}
\max_{s\in{\mathcal
S}} \biggl(\frac{\tr(V_sV_{\xi}^{-1})-\tr(Q_sQ_{\xi}^{-1})}{t-1}+\frac
{t-2}{t-1}
\frac{\ell'V_s\ell}{\ell'V_{\xi}\ell} \biggr)&\leq&1.
\end{eqnarray*}
By taking $\xi_0=\xi$, we have equality in (\ref{eqn8015}). Also
observe that conditioning on fixed~$\xi$, the left-hand side of
(\ref{eqn8015}) is a linear function of the proportions in $\xi_0$.
Hence we have
\begin{eqnarray*}
\max_{s\in{\mathcal
S}} \biggl(\frac{\tr(V_sV_{\xi}^{-1})-\tr(Q_sQ_{\xi}^{-1})}{t-1}+\frac
{t-2}{t-1}
\frac{\ell'V_s\ell}{\ell'V_{\xi}\ell} \biggr)&\geq&1.
\end{eqnarray*}
Then equation (\ref{eqn8016}) follows.
\end{pf}

\begin{remark}\label{remark1}
Theorem \ref{thm1110} neglected the pseudo symmetric measures with
$\det(V_{\xi})=0$ and $\det(Q_{\xi})>0$, and Theorem \ref
{thm1110}(i)--(iii)
neglected those with $\det(Q_{\xi})=0$ and $q_{\xi}(x_{\xi})=0$.
However, all these measures yield $q_{\xi}^*=0$, and thus they cannot be
optimal under ${\mathcal A}_g$, ${\mathcal D}_g$ and ${\mathcal E}_g$
criteria. Actually, $\tau$ is not estimable for any measure with
$q_{\xi
}^*=0$ and
hence such measures should not be adopted [\citet{r28}, Chapter~3]. Note also that $\det(Q_{\xi})=0$ implies $\det(V_{\xi})=0$. Theorem
\ref{thm1110} gives a comprehensive list of conditions to judge the
optimality of a pseudo symmetric measure for estimating $\tau$.
\end{remark}

\begin{remark}
Since we also have $q_{\xi}^*=c_{\xi00}-c_{\xi02}^2/c_{\xi22}$ by
Lemma \ref{lemma1109}, if the function $q_{\xi}(x)$ is replaced with
$c_{\xi00}+2c_{\xi02}x+c_{\xi22}x^2$, equivalent conditions for
optimal pseudo symmetric measures with respect to the four criteria
could be
derived similarly.
\end{remark}

\begin{remark}
For the nonproportional model (\ref{eqn729}), the information matrix
of a pseudo symmetric measure has $t-1$ eigenvalues of $q_{\xi}^*$ and
one of 0. The measure in Theorem \ref{thm1110}(iii) should also be
universally optimal under model (\ref{eqn729}).
\end{remark}

\begin{remark}
When there is only one neighbor effect, say left, we have $\det(Q_{\xi
})=0$ for all $\xi\in{\mathcal P}$. Theorem \ref{thm1110} reduces to
equivalent conditions for the optimal crossover measures where the
pre-period treatment is equal to the treatment in the last period for each
subject.
\end{remark}

\section{Optimal designs for total treatment effect}\label{secOptimaldesignsfortotaleffect}

In this section we study optimal measures for estimating the total
treatment effect, defined by $\theta=(1+\lambda_1+\lambda_2)\tau$.
\citet{r4} commented that the total treatment effect is more
important when the experiment is aimed at finding a single treatment
which is
recommended for use in the whole field.

When $1+\lambda_1+\lambda_2=0$, $\theta$ takes the value of constant
$0$ regardless the value of~$\tau$, and there is no need to carryout the
experiment. In the following we assume $1+\lambda_1+\lambda_2\neq0$.
By plugging $\tau=\theta/(1+\lambda_1+\lambda_2)$ into model
(\ref{eqn728}), we have
\begin{eqnarray*}
\label{eqn7302} Y_d&=&1_{nk}\mu+U\beta+(1+
\lambda_1+\lambda_2)^{-1}(T_d+
\lambda _1L_d+\lambda_2R_d)
\theta+\varepsilon.
\end{eqnarray*}
The information matrix for $\theta$ is
\begin{eqnarray*}
C_d(\theta)&=&(1+\lambda_1+\lambda_2)^{-2}
\\
&&{}\times M_{1,\lambda_1,\lambda_2}'\operatorname{pr}^{\bot} (
\tilde{U}|M_{-1,1+\lambda_2,-\lambda_2}\theta|M_{-1,-\lambda
_1,1+\lambda_1}\theta)M_{1,\lambda_1,\lambda_2}.
\end{eqnarray*}
Here we used the equation $\operatorname{pr}^{\bot}EF=\operatorname{pr}^{\bot
}E$ for any nonsingular matrix
$F$. Actually, it is seen that $1+\lambda_1+\lambda_2=0$ will yield
infinite $C_d(\theta)$ for any $d$, which implies that the covariance
matrix for $\theta$ is zero. Our previous comment on this special case is
justified here. In the same way that we defined $C_{\xi}(\tau)$ in
Section~\ref{secOptimaldesignsfordirecteffect}, the information
matrix of a
measure $\xi$ for $\theta$ is given by $C_{\xi}(\theta
)=n^{-1}C_d(\theta
)$, which is independent of $n$ and can be expressed in a similar
fashion to equation (\ref{eqn2014805}). In the spirit of Proposition
\ref{prop730}, we shall restrict our considerations to pseudo symmetric
measures.

To precede, we define $\ell_0=(-1,1+\lambda_2,-\lambda_2)'$, $\ell
_1=(-1,-\lambda_1,1+\lambda_1)'$, $L_0=(\ell_0,\ell_1)$ and
$L_1=(\ell,\ell_0,\ell_1)$. Let $V_{\xi,1}=L_1'V_{\xi}L_1$,
$Q_{\xi
,1}=L_0'V_{\xi}L_0$ and $q_{\xi,1}^*=\ell'V_{\xi}\ell-\ell'V_{\xi}L_0
Q_{\xi,1}^-L_0'V_{\xi}\ell$.

\begin{proposition}\label{prop04031}
For a pseudo symmetric measure $\xi$, the information matrix $C_{\xi
}(\theta)$ has eigenvalues of $0$,
$(1+\lambda_1+\lambda_2)^{-2}(t-1)^{-1}q_{\xi,1}^*$ and $(1+\lambda
_1+\lambda_2)^{-2}(t-1)^{-1}\ell'V_{\xi}\ell$ with multiplicities of
$1$, $1$
and $t-2$, respectively. Moreover we have $q^*_{\xi,1}\leq\ell
'V_{\xi
}\ell$.
\end{proposition}

\begin{pf}
Denote $\tilde{A}_{d}=(M_{-1,1+\lambda_2,-\lambda_2}\theta
|M_{-1,-\lambda_1,1+\lambda_1}\theta)$. Using $1_t'\theta=0$ and
$C_{\xi ij}=c_{\xi
ij}B_t/(t-1)$, we have
\begin{eqnarray*}
C_{\xi}(\theta)&=& n^{-1}(1+\lambda_1+
\lambda_2)^{-2}
\\
&&{}\times\bigl\{ M_{1,\lambda
_1,\lambda_2}'
\operatorname{pr}^{\bot}(\tilde{U})M_{1,\lambda_1,\lambda
_2}\\
&&\hspace*{17pt}{}-M_{1,\lambda_1,\lambda_2}'\operatorname{pr}^{\bot}(\tilde{U}) \tilde
{A}_{d} \bigl[\tilde {A}_{d}'
\operatorname{pr}^{\bot}(\tilde{U}) \tilde{A}_{d}\bigr]^-
\tilde{A}_{d}'\operatorname{pr}^{\bot}(\tilde{U})
M_{1,\lambda_1,\lambda
_2} \bigr\}
\\
&=&(1+\lambda_1+\lambda_2)^{-2}(t-1)^{-1}
\bigl[\ell'V_{\xi}\ell B_t-a\bigl(\theta
'\theta\bigr)^{-1}\theta\theta'\bigr],
\end{eqnarray*}
where $a=\ell'V_{\xi}L_0 Q_{\xi,1}^-L_0'V_{\xi}\ell$. Let $\{
x_1,\ldots
,x_{t-2}\}$ be the orthogonal basis which is orthogonal to both $1_t$ and
$\theta$. Then $\{1_t,\theta,x_1,\ldots,x_{t-2}\}$ forms the
eigenvectors of $C_{\xi}(\theta)$. The corresponding eigenvalues are $0$,
$(1+\lambda_1+\lambda_2)^{-2}(t-1)^{-1}q_{\xi,1}^*$ and $(1+\lambda
_1+\lambda_2)^{-2}(t-1)^{-1}\ell'V_{\xi}\ell$ with multiplicities of
$1$, $1$
and $t-2$, respectively. The proof is completed in view of $a\ge0$.
\end{pf}

Since $V_{\xi}=0$ implies $V_{\xi,1}=0$, we neglect the measures with
each supporting sequence consisting of only one treatment. Note that
$q_{\xi,1}^*$ is the same Schur complement of $V_{\xi,1}$ as $q_{\xi
}^*$ is that of $V_{\xi}$. Define $V_{s,1}=L_1'V_sL_1$ and
$Q_{s,1}=L_0'V_sL_0$. Let $q_{\xi,1}(x)$ be the same function of
$V_{\xi
,1}$ as $q_{\xi}(x)$ is that of $V_{\xi}$, and $q_{s,1}(x)$ be the same
function of $V_{s,1}$ as $q_{s}(x)$ is that of $V_{s}$. Similar
arguments for Theorem \ref{thm1110} yield the following theorem.

\begin{theorem}\label{thm11102}
In estimating $\theta$ under model (\ref{eqn728}), a pseudo symmetric
measure $\xi$ with $\ell_0'V_{\xi}\ell_0>0$ is optimal in the following
cases. In each case, the $p_s$ in $\xi$ is positive only if $s$ reaches
the maximum therein.
\begin{longlist}[(iii)]
\item[(i)] If $\det(Q_{\xi,1})=0$, then $\xi$ is ${\mathcal A}_g$-optimal if
and only if
\begin{eqnarray*}
\max_{s\in{\mathcal S}}\frac{q_{\xi,1}(x_{\xi,1})^{-2}
q_{s,1}(x_{\xi
,1})+(t-2)(\ell'V_{\xi}\ell)^{-2}\ell'V_s\ell}{
q_{\xi,1}(x_{\xi,1})^{-1}+(t-2)(\ell'V_{\xi}\ell)^{-1}}&=&1.
\end{eqnarray*}
If $\det(V_{\xi,1})>0$, then $\xi$ is ${\mathcal A}_g$-optimal if and
only if
\begin{eqnarray*}
\max_{s\in{\mathcal
S}}\frac{r_{s,1}\det(Q_{\xi,1})/\det(V_{\xi,1})+(t-2)(\ell'V_{\xi
}\ell
)^{-2}\ell'V_s\ell}{\det(Q_{\xi,1})/\det(V_{\xi,1})+(t-2)(\ell
'V_{\xi
}\ell)^{-1}}&=&1.
\end{eqnarray*}
\item[(ii)] If $\det(Q_{\xi,1})=0$, then $\xi$ is ${\mathcal D}_g$-optimal if
and only if
\begin{eqnarray*}
\max_{s\in{\mathcal
S}} \biggl(\frac{1}{t-1}\frac{q_{s,1}(x_{\xi,1})}{q_{\xi,1}(x_{\xi
,1})}+
\frac{t-2}{t-1}\frac{\ell'V_s\ell}{\ell'V_{\xi}\ell} \biggr)&=&1.
\end{eqnarray*}
If $\det(V_{\xi,1})>0$, then $\xi$ is ${\mathcal D}_g$-optimal if and
only if
\begin{eqnarray*}
\max_{s\in{\mathcal S}} \biggl(\frac{r_{s,1}}{t-1}+\frac
{t-2}{t-1}
\frac
{\ell'V_s\ell}{\ell'V_{\xi}\ell} \biggr)&=&1.
\end{eqnarray*}
\item[(iii)] If $\det(Q_{\xi,1})=0$, then $\xi$ is ${\mathcal E}_g$-optimal if
and only if
\begin{eqnarray*}
\max_{s\in{\mathcal S}}\frac{q_{s,1}(x_{\xi,1})}{q_{\xi,1}(x_{\xi
,1})}&=&1.
\end{eqnarray*}
If $\det(V_{\xi,1})>0$, then $\xi$ is ${\mathcal E}_g$-optimal if and
only if
\begin{eqnarray*}
\max_{s\in{\mathcal S}}r_{s,1}&=&1.
\end{eqnarray*}
\item[(iv)] If $\det(Q_{\xi,1})=0$, then $\xi$ is ${\mathcal T}_g$-optimal if
and only if
\begin{eqnarray*}
\max_{s\in{\mathcal S}}\frac{q_{s,1}(x_{\xi,1})+(t-2)\ell'V_s\ell
}{q_{\xi,1}(x_{\xi,1})+(t-2)\ell'V_{\xi}\ell}&=&1.
\end{eqnarray*}
If $\det(V_{\xi,1})>0$, then $\xi$ is ${\mathcal E}_g$-optimal if and
only if
\begin{eqnarray*}
\max_{s\in{\mathcal S}}\frac{r_{s,1} \det(V_{\xi,1})/\det(Q_{\xi
,1})+(t-2)\ell'V_s\ell}{\det(V_{\xi,1})/\det(Q_{\xi,1})+(t-2)\ell
'V_{\xi
}\ell}&=&1.
\end{eqnarray*}
Here $r_{s,1}=\tr(V_{s,1}V_{\xi,1}^{-1})-\tr(Q_{s,1}Q_{\xi,1}^{-1})$ and
$x_{\xi,1}=-\ell'V_{\xi}\ell_0/\ell_0'V_{\xi}\ell_0$.
\end{longlist}
\end{theorem}

\begin{remark}
With arguments similar to those in Remark \ref{remark1}, Theorem \ref
{thm11102} gives a comprehensive list of conditions to judge the
optimality of
pseudo symmetric measures with $\ell_0' V_{\xi}\ell_0>0$ for estimating
$\theta$. Equivalence conditions for the four criteria could be easily
derived when $\ell_0' V_{\xi}\ell_0=0$, where we need to consider
whether the cases of $\ell_1' V_{\xi}\ell_1$ are equal to 0 or not,
separately. We omit
the details due to limit of space.
\end{remark}

\begin{remark}
If the within-block covariance matrix $\Sigma$ is of type-H, $\ell_0'
V_{\xi}\ell_0=0$ implies $V_{\xi}\ell_0=0$, and thus $V_{\xi}=0$ in
view of
equation (\ref{Vs}) below. Similarly, $\ell_1' V_{\xi}\ell_1=0$ also
results in $V_{\xi}=0$. Therefore, except for measures with each
supporting sequence consisting of only one treatment, we have $\ell_0'
V_{\xi}\ell_0>0$ and $\ell_1' V_{\xi}\ell_1>0$ for any type-H
matrix $\Sigma$.
\end{remark}

The following proposition shows that the values of $\lambda_1$ and
$\lambda_2$ are irrelevant to the determination of the
$\mathcal{E}_g$-optimal pseudo symmetric measures for estimating
$\theta$.

\begin{proposition}\label{Proposition0613}
For any measure, we have
\begin{eqnarray*}
&& q_{\xi,1}^*=(1+\lambda_1+\lambda_2)^2
\min_{x,y}\bigl[(1,x,y)\Gamma' V_{\xi
}
\Gamma(1,x,y)'\bigr],
\end{eqnarray*}
where
\begin{eqnarray*}
\Gamma=\pmatrix{
1 & -1 & -1
\cr
0 & 1 & 0
\cr
0 & 0 & 1}.
\end{eqnarray*}
\end{proposition}
\begin{pf}
Let
\begin{eqnarray*}
\Lambda=(1+\lambda_1+\lambda_2)^{-1}\pmatrix{
1+\lambda_1+\lambda_2 & 0
& 0
\cr
\lambda_1 & 1+\lambda_2 & -\lambda_1
\cr
\lambda_2 & -\lambda_2 & 1+\lambda_1}.
\end{eqnarray*}
Note that $L_1=(\ell,\ell_0,\ell_1)=(1+\lambda_1+\lambda_2)\Gamma
\Lambda$. From Proposition 3 in \citet{r24}, we have
\begin{eqnarray*}
q_{\xi,1}^*&=&\min_{x,y}\bigl[(1,x,y)
V_{\xi,1}(1,x,y)'\bigr]
\\
&=&(1+\lambda_1+\lambda_2)^2\min
_{x,y}\bigl[(1,x,y)\Lambda'\Gamma'
V_{\xi
}\Gamma\Lambda(1,x,y)'\bigr]
\\
&=&(1+\lambda_1+\lambda_2)^2\min
_{x,y}\bigl[(1,x,y)\Gamma' V_{\xi
}
\Gamma(1,x,y)'\bigr].
\end{eqnarray*}
The last equality uses the fact that for all possible values of $x$ and
$y$, $\Lambda(1,x,y)'$ and $(1,x,y)'$ share the same vector space.
\end{pf}

\section{Optimal designs for the undirectional model}\label
{secUndirectionalsideeffects}

In many applications, it is reasonable to assume $\lambda_1=\lambda
_2:=\lambda$; that is, the neighbor effects do not depend on whether
they are from
the left or right. See \citet{r11}, \citet{r7} and \citet{r15}, for examples. Under this condition, model
(\ref{eqn728}) reduces to
\begin{eqnarray}
\label{eqn730} Y_d&=&1_{nk}\mu+U\beta+(T_d+
\lambda L_d+\lambda R_d)\tau+\varepsilon.
\end{eqnarray}
The information matrix for $\tau$ under model (\ref{eqn730}) is
\begin{eqnarray*}
\tilde C_{d}(\tau)&=&M_{1,\lambda,\lambda}'
\operatorname{pr}^{\bot
}(\tilde {U}|M_{0,1,1}\tau)M_{1,\lambda,\lambda}.
\end{eqnarray*}
The information matrix of a measure $\xi$ for $\tau$ is $\tilde
C_{\xi
}(\tau)=n^{-1}\tilde C_{d}(\tau)$. Also we consider only optimal measures
in the pseudo symmetric format.

Define $\ell_2=(0,1,1)'$ and $L_2=(\ell,\ell_2)$, where $\ell$ is
defined in Section~\ref{secOptimaldesignsfordirecteffect} with the
value of
$(1,\lambda,\lambda)'$ here. Let $V_{\xi,2}=L_2'V_{\xi}L_2$,
$Q_{\xi,2}=\ell_2'V_{\xi}\ell_2$, $V_{s,2}=L_2'V_{s}L_2$,
$Q_{s,2}=\ell
_2'V_{s}\ell_2$ and
$q_{\xi,2}^*=\ell'V_{\xi}\ell-\ell'V_{\xi}\ell_2Q_{\xi,2}^-\ell
_2'V_{\xi
}\ell$. Similar to Proposition \ref{prop04031}, we have the following.

\begin{proposition}\label{prop0403}
For a pseudo symmetric measure $\xi$, the information matrix $\tilde
C_{\xi}(\tau)$ has eigenvalues of $0$, $(t-1)^{-1}q_{\xi,2}^*$ and
$(t-1)^{-1}\ell'V_{\xi}\ell$ with multiplicities of $1$, $1$ and $t-2$,
respectively. Moreover we have $q_{\xi,2}^*\le\ell'V_{\xi}\ell$.
\end{proposition}

Note that if $Q_{\xi,2}=0$, then $q_{\xi,2}^*=\ell'V_{\xi}\ell
=c_{\xi
00}$ and hence $\tilde C_{\xi}(\tau)=c_{\xi00}B_t/(t-1)$. By arguments
similar to those in Theorem \ref{thm1110}, we obtain the following result.

\begin{theorem}\label{thm20140801}
In estimating $\tau$ under model (\ref{eqn730}), a pseudo symmetric
measure $\xi$ is optimal in the following cases. In each case, the
$p_s$ in
$\xi$ is positive only if $s$ reaches the maximum therein.
\begin{longlist}[(ii)]
\item[(i)] If $Q_{\xi,2}=0$, then $\xi$ is universally optimal if and only if
\begin{eqnarray*}
\max_{s\in{\mathcal S}}\frac{c_{s 00}}{c_{\xi00}}&=&1.
\end{eqnarray*}
\item[(ii)] If $\det(V_{\xi,2})>0$, then $\xi$ is ${\mathcal A}_g$-optimal if
and only if
\begin{eqnarray*}
\max_{s\in{\mathcal S}}\frac{r_{s,2} Q_{\xi,2}/\det(V_{\xi
,2})+(t-2)(\ell'V_{\xi}\ell)^{-2}\ell'V_s\ell}{ Q_{\xi,2}/\det(
V_{\xi,2})+(t-2)(\ell'V_{\xi}\ell)^{-1}}&=&1.
\end{eqnarray*}
$\xi$ is ${\mathcal D}_g$-optimal if and only if
\begin{eqnarray*}
\max_{s\in{\mathcal S}} \biggl(\frac{r_{s,2}}{t-1}+\frac
{t-2}{t-1}
\frac
{\ell'V_s\ell}{\ell'V_{\xi}\ell} \biggr)&=&1.
\end{eqnarray*}
$\xi$ is ${\mathcal E}_g$-optimal if and only if
\begin{eqnarray*}
\max_{s\in{\mathcal S}} r_{s,2}&=&1.
\end{eqnarray*}
$\xi$ is ${\mathcal T}_g$-optimal if and only if
\begin{eqnarray*}
\max_{s\in{\mathcal S}}\frac{ r_{s,2} \det( V_{\xi,2})/ Q_{\xi
,2}+(t-2)\ell'V_s\ell}{\det( V_{\xi,2})/ Q_{ \xi,2}+(t-2)\ell
'V_{\xi
}\ell}&=&1.
\end{eqnarray*}
\item[(iii)] Otherwise, $\xi$ is not optimal.

Here $r_{s,2}=\tr(V_{s,2} V_{\xi,2}^{-1})- Q_{s,2}Q_{\xi,2}^{-1}$.
\end{longlist}
\end{theorem}

It is easy to verify that for any measure,
\begin{eqnarray}
\label{eqn0613} q_{\xi,2}^*&=&\min_{x}
\bigl[(1,x)V_{\xi,2}(1,x)'\bigr]=\min_{x}
\bigl[(1,x,x)V_{\xi
}(1,x,x)'\bigr].
\end{eqnarray}
Therefore, the value of $\lambda$ is irrelevant to the search of
${\mathcal E}_g$-optimal pseudo symmetric measures for estimating $\tau
$.

Next, we consider the total treatment effect $\theta$. With the reason
we explained earlier, we shall assume $1+2\lambda\neq0$. The information
matrix for $\theta$ under model~(\ref{eqn730}) is
\begin{eqnarray*}
\tilde C_{d}(\theta)&=&(1+2\lambda)^{-2}M_{1,\lambda,\lambda
}'
\operatorname{pr}^{\bot} (\tilde{U}|M_{2,-1,-1}\theta)M_{1,\lambda,\lambda}.
\end{eqnarray*}
For a measure $\xi$, its information matrix for $\theta$ is given by
$\tilde C_{\xi}(\theta)=n^{-1}\tilde C_{d}(\theta)$. Now we define
$\ell_3=(2,-1,-1)'$ and $L_3=(\ell,\ell_3)$. Let $V_{\xi
,3}=L_3'V_{\xi
}L_3$, $Q_{\xi,3}=\ell_3'V_{\xi}\ell_3$, $V_{s,3}=L_3'V_{s}L_3$,
$Q_{s,3}=\ell_3'V_{s}\ell_3$ and $q_{\xi,3}^*=\ell'V_{\xi}\ell
-\ell
'V_{\xi}\ell_3Q_{\xi,3}^-\ell_3'V_{\xi}\ell$. Also we have the following.

\begin{proposition}\label{prop0610}
For a pseudo symmetric measure $\xi$, the information matrix $\tilde
C_{\xi}(\theta)$ has eigenvalues of $0$,
$(1+2\lambda)^{-2}(t-1)^{-1}q_{\xi,3}^*$ and $(1+2\lambda
)^{-2}(t-1)^{-1}\ell'V_{\xi}\ell$ with multiplicities of $1$, $1$
and $t-2$,
respectively. Moreover we have $q_{\xi,3}^*\le\ell'V_{\xi}\ell$.
\end{proposition}

Note that if $Q_{\xi,3}=0$, then $q_{\xi,3}^*=\ell'V_{\xi}\ell
=(1+2\lambda)^2c_{\xi00}$ and hence $\tilde C_{\xi}(\theta)=c_{\xi
00}B_t/(t-1)$.
Similar to Theorem \ref{thm20140801}, we obtain the following theorem.

\begin{theorem}\label{thm20140802}
In estimating $\theta$ under model (\ref{eqn730}), a pseudo symmetric
measure $\xi$ is optimal in the following cases. In each case, the $p_s$
in $\xi$ is positive only if $s$ reaches the maximum therein.
\begin{longlist}[(iii)]
\item[(i)]  If $Q_{\xi,3}=0$, then $\xi$ is universally optimal if and only if
\begin{eqnarray*}
\max_{s\in{\mathcal S}}\frac{c_{s 00}}{c_{\xi00}}&=&1.
\end{eqnarray*}
\item[(ii)]  If $\det(V_{\xi,3})>0$, then $\xi$ is ${\mathcal A}_g$-optimal if
and only if
\begin{eqnarray*}
\max_{s\in{\mathcal S}}\frac{r_{s,3} Q_{\xi,3}/\det(V_{\xi
,3})+(t-2)(\ell'V_{\xi}\ell)^{-2}\ell'V_s\ell}{ Q_{\xi,3}/\det(
V_{\xi,3})+(t-2)(\ell'V_{\xi}\ell)^{-1}}&=&1.
\end{eqnarray*}
$\xi$ is ${\mathcal D}_g$-optimal if and only if
\begin{eqnarray*}
\max_{s\in{\mathcal S}} \biggl(\frac{r_{s,3}}{t-1}+\frac
{t-2}{t-1}
\frac
{\ell'V_s\ell}{\ell'V_{\xi}\ell} \biggr)&=&1.
\end{eqnarray*}
$\xi$ is ${\mathcal E}_g$-optimal if and only if
\begin{eqnarray*}
\max_{s\in{\mathcal S}} r_{s,3}&=&1.
\end{eqnarray*}
$\xi$ is ${\mathcal T}_g$-optimal if and only if
\begin{eqnarray*}
\max_{s\in{\mathcal S}}\frac{ r_{s,3} \det( V_{\xi,3})/ Q_{\xi
,3}+(t-2)\ell'V_s\ell}{\det( V_{\xi,3})/ Q_{ \xi,3}+(t-2)\ell
'V_{\xi
}\ell}&=&1.
\end{eqnarray*}
\item[(iii)]  Otherwise, $\xi$ is not optimal.

Here $r_{s,3}=\tr(V_{s,3} V_{\xi,3}^{-1})- Q_{s,3}Q_{\xi,3}^{-1}$.
\end{longlist}
\end{theorem}

It is easy to verify that for any measure,
\begin{eqnarray}
\label{eqn06132} &&\quad q_{\xi,3}^*=\min_{x}
\bigl[(1,x)V_{\xi,3} (1,x)'\bigr]=(1+2\lambda)^2
\min_{x}\bigl[(1,x,x)\Gamma'V_{\xi}
\Gamma(1,x,x)'\bigr].
\end{eqnarray}
Therefore, the value of $\lambda$ is also irrelevant in the search of
${\mathcal E}_g$-optimal pseudo symmetric measures for estimating
$\theta$.

Finally, we establish the connection between optimal measures for the
directional and undirectional models if the within-block covariance matrix
$\Sigma$ is of \mbox{type-H}.

\begin{lemma}\label{lemma0807}
If $\Sigma$ is of type-H, we have $q_{\xi,2}^*=q_{\xi}^*$ and
$(1+2\lambda)^{-2}q_{\xi,3}^*=(1+\lambda_1+\lambda_2)^{-2}q_{\xi,1}^*$.
\end{lemma}
\begin{pf}
Note that $\tilde{B}=B_k$ if $\Sigma$ is of type-H. For a sequence
$s=(t_1,\ldots,t_k)$, define $t_0=t_k$ and $t_{k+1}=t_1$. Let $k_j$ be the
frequency of treatment $i$ appearing in $s$. Clearly, $\sum_{i=1}^t
k_i=k$. Let $m_s=k^{-1}\sum_{i=1}^t k_i^2$, $f_{s}=\sum_{i=1}^k
\mathbb{I}_{t_{i}=t_{i-1}}$, $g_{s}=\sum_{i=1}^k \mathbb
{I}_{t_{i}=t_{i+1}}$ and
$h_{s}=\sum_{i=1}^k \mathbb{I}_{t_{i-1}=t_{i+1}}$.
By straightforward calculations, we have
$c_{s00}=c_{s11}=c_{s22}=k-m_s$, $c_{s01}=f_s-m_s$, $c_{s02}=g_s-m_s$
and $c_{s12}=h_s-m_s$.
Since $f_s=g_s$, we have
\begin{eqnarray}
\label{Vs}
V_s=(c_{sij})_{0\leq i,j\leq2}=\pmatrix{
k-m_s & f_s-m_s
& f_s-m_s
\cr
f_s-m_s & k-m_s & h_s-m_s
\cr
f_s-m_s & h_s-m_s &
k-m_s}.
\end{eqnarray}
Note that $V_s=0$ if and only if $s$ consists of only one treatment.
From Proposition~3 in \citet{r24}, we have
$q_{\xi}^*=\min_{x,y}[(1,x,y)V_{\xi}(1,x,y)']$. Since $(1,x,y)V_{\xi
}(1,x,y)'$ is convex and exchangeable in $x$ and $y$ by equation (\ref
{Vs}) and
$V_{\xi}=\sum_{s\in\mathcal S}p_s V_s$, it can
achieve the minimum at some point of $x=y$. Therefore, $q_{\xi}^*=\min_{x}(1,x,x)V_{\xi}(1,x,x)'=q_{\xi,2}^*$ in view of equation (\ref{eqn0613}).

From Proposition \ref{Proposition0613}, $q_{\xi,1}^*=(1+\lambda
_1+\lambda_2)^2\min_{x,y}[(1,x,y)\Gamma' V_{\xi}\Gamma(1,x,y)']$.
By equation (\ref{Vs}) and
$V_{\xi}=\sum_{s\in\mathcal S}p_s V_s$, we know $(1,x,y)\Gamma'
V_{\xi
}\Gamma(1,x,y)'$ is convex and exchangeable in $x$ and $y$. Thus it can
achieve the minimum at some point of $x=y$. Then $(1+\lambda_1+\lambda
_2)^{-2}q_{\xi,1}^*=\min_{x}[(1,x,x)\Gamma' V_{\xi}\Gamma
(1,x,x)']=(1+2\lambda)^{-2}q_{\xi,3}^*$
%$(1+2\lambda)^{-2}q_{\xi,3}^*=(1+\lambda_1+\lambda_2)^{-2}q_{\xi,1}^*$
in view of equation (\ref{eqn06132}).
\end{pf}

\begin{theorem}\label{thm20140803}
If $\Sigma$ is of type-H, a pseudo symmetric measure is ${\mathcal
E}_g$-optimal for $\tau$ (resp., $\theta$) under model (\ref{eqn730})
if and only
if it is ${\mathcal E}_g$-optimal for $\tau$ (resp., $\theta$) under
model (\ref{eqn728}). Furthermore, if $\lambda_1=\lambda_2=\lambda$,
the same
result holds for ${\mathcal A}_g$-, ${\mathcal D}_g$- and ${\mathcal
T}_g$-optimal pseudo symmetric measures.
\end{theorem}

This theorem is readily proved by using Lemma \ref{lemma0807} and
Propositions \ref{prop7302}, \ref{prop04031}, \ref{prop0403}~and~\ref{prop0610}.

\section{Examples}\label{secexamples}
For a sequence $s=(t_1,\ldots,t_k)$, define the symmetric block of~$s$
as $\langle s \rangle=\{\sigma s\dvtx  \sigma\in{\mathcal O}\}$. A symmetric
block is an equivalence class, and hence ${\mathcal S}$ is partitioned
into $m+1$ symmetric blocks, say $\langle s_0\rangle,\langle
s_1\rangle,\ldots,\langle s_m\rangle$, where $s_i$'s are the
representative sequences in their own blocks. Without loss of
generality, let
$\langle s_0\rangle$ be the symmetric block of sequences with identical
elements. For a measure $\xi=(p_s,s\in\mathcal S)$, let $p_{\langle
s_i\rangle}=\sum_{s\in\langle s_i \rangle}p_{s}$ and $P_{\xi
}=(p_{\langle s_1 \rangle},\ldots,p_{\langle s_m \rangle})$. Since
$V_{s}$ is
invariant for sequences in the same symmetric block, two pseudo
symmetric measures with the same $P_{\xi}$ will share the same value of
$\phi_{g}(\xi)$. By Remark 2 in \citet{r31}, one can derive an exact
optimal design in two steps: First, find the optimal $P_{\xi}$, and then
construct an exact pseudo symmetric design with that $P_{\xi}$ by using
some combinatory structures, such as type I orthogonal arrays [\citet{r29}]. See \citet{r3} and \citet{r5} for more techniques to construct exact pseudo symmetric designs.

Note that $C_{sij}=0$, $0\le i,j\le2$, for any $s\in\langle
s_0\rangle
$. Given a measure $\xi$ with $p_{\langle s_0\rangle}>0$, one can always
obtain a measure superior to $\xi$ by replacing all sequences in
$\langle s_0\rangle$ with sequences not in the set.
%and raising proportions of other sequences at the ratio of
%$(1-p_s)^{-1}$.
Therefore, the symmetric block $\langle s_0\rangle$ will be ignored in
the following discusssion.

In the sequel, we will determine the optimal $P_{\xi}$ under model
(\ref
{eqn728}) through computer search based on Theorems \ref{thm1110} and
\ref{thm11102}. The one for the undirectional model~(\ref{eqn730})
can be determined in a similar way by using Theorems \ref{thm20140801},
\ref{thm20140802} and \ref{thm20140803}. The general algorithm for
deriving the optimal $P_{\xi}$ can be obtained by small modifications of
the algorithm in \citet{r32}. For ease of illustration, we consider
only $2\le t, k\le5$ and use the within-block covariance matrix to be of
the form $\Sigma=(\mathbb{I}_{i=j}+\rho\mathbb{I}_{i-j=\pm1
(\operatorname{mod} \, \,
k)})_{1\le i,j\le k}$. In the following examples, we take $\rho$ in
$\{0,-0.3, 0.3\}$. Note that $\rho=0$ implies $\Sigma=I_k$; that is,
the errors are uncorrelated. First, let $\lambda_1$ and $\lambda_2$ be
nonnegative values from $[0,1]$, and the negative case will be
discussed later. All measures given below are pseudo symmetric measures.
\begin{longlist}[{}]
\item[\textit{Cases} $k=2$ and $3$.] When $k=2$, the symmetric block is
$\langle
12 \rangle$. When $k=3$, the symmetric block is $\langle112
\rangle$ for $t=2$, and those are $\langle112 \rangle$ and $\langle
123 \rangle$ for $t\ge3$. By straightforward calculations, it can be
verified that the second smallest eigenvalues of $C_{\xi}(\tau)$ and
$C_{\xi}(\theta)$ are both zero for any measure when $k=2$ and $3$.
Therefore, neither $\tau$ nor $\theta$ is estimable, and the optimal
measures do not exist in these cases. This phenomenon is also observed by
\citet{r4} and \citet{r14} for the
nonproportional interference model.

\item[\textit{Case of} $k=4$.] When $t=2$,
%since the information matrix of $\tau$ or $\theta$ has only one
%nonzero eigenvalues,
the four criteria become the same one. From Propositions \ref
{prop7302}, \ref{prop04031} and \ref{Proposition0613}, it is known
that the
optimality of a measure for $\tau$ or $\theta$ does not depend on the
values of $\lambda_1$ and $\lambda_2$. For $\tau$, we find that the measure
with $p_{\langle1122 \rangle}=1$ is optimal for the three values of
$\rho$. Next, consider the optimal measures for $\theta$. If $\rho
=0$, the
measure with $p_{\langle1122 \rangle}=2/3$ and $p_{\langle1212
\rangle
}=1/3$ is optimal. An exact pseudo symmetric design with three blocks
based on it is given by
\begin{eqnarray*}
\pmatrix{ 1 & 1 & 2 & 2
\cr
1 & 1 & 2 & 2
\cr
1 & 2 & 1 & 2 }.
\end{eqnarray*}
In order to implement this design in practice, one has to adopt the
randomization procedure as suggested by \citet{r3}. If $\rho=-0.3$, the measure with $p_{\langle1122 \rangle
}=0.61$ and $p_{\langle1212 \rangle}=0.39$ is optimal. If $\rho=0.3$,
the measure with
$p_{\langle1122 \rangle}=0.76$ and $p_{\langle1212 \rangle}=0.24$
is optimal.

When $t=3$, the measure with $p_{\langle1123 \rangle}=1$ is optimal
for $\tau$ under the four criteria given all the values of $\lambda_1$,
$\lambda_2$ and $\rho$. Consider the optimal measures for $\theta$. If
$\rho=0$, the ${\mathcal A}_g$-, ${\mathcal D}_g$- and ${\mathcal
T}_g$-optimal
measures vary for different values of $\lambda_1$ and $\lambda_2$. For
all of them, there are two supporting symmetric blocks, that is,
$\langle
1123 \rangle$ and $\langle1213 \rangle$. Meanwhile, the former
symmetric block dominates. The measure with $p_{\langle1123 \rangle
}=2/3$ and
$p_{\langle1213 \rangle}=1/3$ is $\mathcal{E}_g$-optimal. If $\rho
=-0.3$ and $0.3$, we observe that the supporting symmetric blocks are the
same as those for $\rho=0$, except for the proportions of the
supporting symmetric blocks.

When $t=4$ and $5$, we find that the measure with $p_{\langle1234
\rangle}=1$ is optimal for both $\tau$ and $\theta$ under the four criteria,
given all the values of $\lambda_1$, $\lambda_2$ and $\rho$.

\item[\textit{Case of} $k=5$.]  When $t=2$, for both $\tau$ and $\theta$ we have
the following. The measure with $p_{\langle11122 \rangle}=0.8$ and
$p_{\langle11212 \rangle}=0.2$ is optimal for $\rho=0$,
the measure with $p_{\langle11122 \rangle}=0.71$ and $p_{\langle
11212 \rangle}=0.29$ is optimal for $\rho=-0.3$
and the measure with $p_{\langle11122 \rangle}=0.90$ and $p_{\langle
11212 \rangle}=0.10$ is optimal for $\rho=0.3$.

When $t=3$, first consider optimal measures for $\tau$. If $\rho=0$,
the optimal measures vary for different values of $\lambda_1$ and
$\lambda_2$ while the supporting symmetric blocks are always $\langle
11223 \rangle$ and $\langle12123 \rangle$. The proportion of
$\langle
11223 \rangle$ is almost one for ${\mathcal A}_g$-, ${\mathcal D}_g$-
and ${\mathcal T}_g$-optimal measures and is $0.90$ for the ${\mathcal
E}_g$-optimal
measure. If $\rho=-0.3$, the supporting symmetric blocks remain the
same as those for $\rho=0$ and $\langle11223 \rangle$ still
dominates. If
$\rho=0.3$, the measure with $p_{\langle11223 \rangle}=1$ is optimal
under the four criteria for $\lambda_1,\lambda_2\in[0,1]$. For
$\theta$,
we have observations similar to those for $\tau$.

When $t=4$, the supporting symmetric blocks are $\langle11234 \rangle$
and $\langle11223 \rangle$. When $t=5$, the supporting symmetric blocks
are $\langle11234 \rangle$, $\langle11223 \rangle$ and $\langle12345
\rangle$. The optimal proportions and the dominating block may change
for different values of $\lambda_1$, $\lambda_2$ and $\rho$.

From Theorems \ref{thm1110} and \ref{thm11102}, it is seen that as
$t$ increases, the equivalent conditions for optimal measures under
${\mathcal
A}_g$, ${\mathcal D}_g$ and ${\mathcal T}_g$ criteria tend to agree
with each other. For example, take $k=5$, $\rho=0$, $\lambda_1=0.1$ and
$\lambda_2=0.2$. The measure with $p_{\langle12345 \rangle}=1$ is
optimal under the three criteria for both $\tau$ and $\theta$ when
$t\ge12$.
Meanwhile, the measure with $p_{\langle11234 \rangle}=0.955$ and
$p_{\langle12345 \rangle}=0.045$ is $\mathcal{E}_g$-optimal for both
$\tau$
and $\theta$ as long as $t\ge5$.

Though the four criteria do not lead to the same optimal measure in
general, the optimal measure under one criterion is typically highly efficient
under the other three. Here the efficiency of a measure under a
criterion is defined as the ratio of $\phi_g(\xi)$ to the maximum
value among
all measures. For the case of $k=5$, $t=3$, $\rho=0$, $\lambda_1=0.1$
and $\lambda_2=0.2$, the efficiencies of optimal measures for $\tau$ are
shown in Table~\ref{table-1} and those for $\theta$ under the four
criteria are shown in Table~\ref{table-2}. They all have efficiencies higher
than 0.97. Furthermore, from the two tables we observe that the optimal
measures for $\tau$ also have high efficiencies in estimating $\theta$
since they are almost the same as those for $\theta$.
\begin{table}[t]
\tablewidth=200pt
\caption{Efficiencies of optimal measures for $\tau$ at
$(k,t,\rho,\lambda_1,\lambda_2)=(5,3,0,0.1,0.2)$}\label{table-1}
\begin{tabular*}{200pt}{@{\extracolsep{\fill}}lcccc@{}}
\hline
$\bolds{p_{\langle11223 \rangle}}$ & \textbf{A} & \textbf{D} & \textbf{E} & \textbf{T} \\
\hline
0.98 & 1\phantom{.00000} & 0.99997 & 0.98817 & 0.99988 \\
0.99 & 0.99998 & 1\phantom{.00000} &0.98670 & 0.99996 \\
0.90 & 0.99265 & 0.99213 & 1\phantom{.00000} & 0.99156\\
1\phantom{.00} & 0.99988 & 0.99997 & 0.98496 & 1\phantom{.00000} \\
\hline
\end{tabular*}
\end{table}

\begin{table}[b]
\tablewidth=200pt
\caption{Efficiencies of optimal measures for $\theta$ at
$(k,t,\rho,\lambda_1,\lambda_2)=(5,3,0,0.1,0.2)$}
\label{table-2}
\begin{tabular*}{200pt}{@{\extracolsep{\fill}}lcccc@{}}
\hline
$\bolds{p_{\langle11223 \rangle}}$ & \textbf{A} & \textbf{D} & \textbf{E} & \textbf{T} \\
\hline
0.93 & 1\phantom{.00000} & 0.99676 & 0.98828 & 0.98702 \\
0.99 & 0.99556 & 1\phantom{.00000} &0.98671 & 0.99787 \\
0.90 & 0.99859 & 0.99213 & 1\phantom{.00000} & 0.97925\\
1\phantom{.00} & 0.99307 & 0.99981 & 0.98215 & 1\phantom{.00000} \\
\hline
\end{tabular*}
\end{table}

From a practical viewpoint, the optimal proportions are sometimes too
harsh for deriving exact designs. However, since the four criterion
functions are continuous in the proportions, we could get a measure
with good proportions in the neighborhood of the optimal one at the
cost of
a little efficiency. For example, when $k=t=5$, $\rho=0$, $\lambda
_1=0.1$ and $\lambda_2=0.2$, the $\mathcal{A}_g$-optimal measure for
$\tau$ is
given by $p_{\langle11223 \rangle}=0.06$ and $p_{\langle12345
\rangle
}=0.94$, which requires $n$ to be a multiple of 50 at least. By rounding
the proportions, we obtain a measure with $p_{\langle12345 \rangle
}=1$, which has efficiency higher than 0.99. An exact pseudo symmetric design
with four blocks based on it is given by \citet{r3} as
\begin{eqnarray*}
&&\pmatrix{ 1 & 2 & 3 & 4 & 5
\cr
1 & 3 & 5 & 2 & 4
\cr
1 & 4 & 2 & 5 & 3
\cr
1 & 5 & 4 & 3 & 2 }.
\end{eqnarray*}

In some occasions, the values of $\lambda_1$ and $\lambda_2$ could be
negative. For example, a good fertilizer will possibly make a plant grow
well so that the plant will compete with its neighbors for the
sunlight, water and other resources in the soil. Suppose
$\lambda_1,\lambda_2\in[-1,0)$. The ${\mathcal A}_g$-, ${\mathcal
D}_g$- and ${\mathcal T}_g$-optimal measures found by the computer
program are different
from the preceding ones for $\lambda_1,\lambda_2\in[0,1]$ in some
cases. In estimating both $\tau$ and $\theta$, we observe the
following. When
$(k,t)=(4,3)$, the supporting symmetric blocks are $\langle1123
\rangle
$ and $\langle1213 \rangle$ for ${\mathcal A}_g$- and ${\mathcal D}_g$-optimal
measures, and are $\langle1212 \rangle$ and $\langle1213 \rangle$ for
${\mathcal T}_g$-optimal measures. Contrarily, for
$\lambda_1,\lambda_2\in[0,1]$, there is only one supporting symmetric
block $\langle1123 \rangle$ for optimal measures in estimating $\tau$,
and the ${\mathcal T}_g$-optimal measure in estimating $\theta$ has two
supporting symmetric blocks as $\langle1123 \rangle$ and $\langle1213
\rangle$. When $(k,t)=(4,4)$, the ${\mathcal A}_g$-optimal measure is
still given by $p_{\langle1234 \rangle}=1$. The supporting symmetric blocks
are $\langle1212 \rangle$, $\langle1213 \rangle$ and $\langle1234
\rangle$ for ${\mathcal D}_g$-optimal measures, and are $\langle1212
\rangle$
and $\langle1213 \rangle$ for ${\mathcal T}_g$-optimal measures. But
for $\lambda_1,\lambda_2\in[0,1]$, there is only one supporting symmetric
block $\langle1234 \rangle$ for optimal measures under the three
criteria. The details for other combinations of parameters are omitted
due to the
limit of space.
\end{longlist}

\section{Discussions}\label{secconclusion}
In this article, two proportional interference models are considered,
in which the neighbor effects of a treatment are proportional to its
direct effect. We investigate the optimal circular designs for the
direct and total treatment effects. Kiefer's equivalence theorems with respect
to A, D, E and T criteria are established, based on which the search of
optimal designs is easy to perform. Moreover, the connection between
optimal designs for the two models is built. Examples are given to
illustrate these theorems for several combinations of parameters.

We now remark on directions for future work. Note that the number of
distinct symmetric blocks will increase, at least geometrically, as the
block size $k$ grows. In such circumstance, it is unlikely that we
could find the optimal proportions within a reasonable amount of time
by using
the current algorithm. Therefore, determining the forms of supporting
symmetric blocks theoretically is vital to solving this problem. As a design
theorist, the ultimate goal is to provide efficient or even optimal
exact designs for any number of blocks. One way to achieve this is to
further explore the constructions of exact pseudo symmetric designs.
The other is to develop methods to build up efficient exact designs by
modifying existing designs of smaller or larger size.

%\section*{Acknowledgements}
%The authors thank the editor, the associate editors, and two referees
%for their valuable comments and suggestions
%that improved the presentation of this article.

\section*{Acknowledgements}
The authors sincerely thank the Editor, the Associate  Editor and two
referees for their insightful comments
and suggestions which have led to the improvement of the paper.

% imsref loaded by daiva.urboniene, 2015-03-02 10:40:06

%\begin{appendix}
%\section{}
%\end{appendix}

%\section*{Acknowledgments}

%\begin{supplement}[id=suppA]
%\slink[doi]{10.1214/00-AOSXXXXSUPP} %[doi,text={...}] - jei reikia
%suskaldyti doi
%\sfilename{aosXXXX\_supp.pdf}
%\end{supplement}

%\begin{thebibliography}{99}
%\bibitem[\protect\citeauthoryear{}{}]{r1}
%\bibitem{r1}
%\end{thebibliography}

\printaddresses
\end{document}